\documentclass[11pt,oneside,a4paper,leqno]{amsart}
\usepackage{amsmath,amsfonts,amsthm}
\usepackage[utf8]{inputenc}
\usepackage{lscape}
\usepackage[T1]{fontenc}
\usepackage{setspace}
\usepackage{tikz}
\usetikzlibrary{positioning}
\usepackage{faktor}
\usepackage{amssymb}
\usepackage{enumitem}
\usepackage{mathtools}
\usepackage{hyperref}
\usepackage{xcolor}
\usepackage[all]{xy}
\usepackage{vmargin}
\usepackage{float}
\usepackage{nomencl}

\setmarginsrb{4.05cm}{3.5cm}{4.05cm}{3cm}{0pt}{1cm}{0pt}{1cm}
\newtheorem{theorem}{Theorem}[section]
\newtheorem*{theorem*}{Theorem}
\newtheorem{corollaire}[theorem]{Corollary}
\newtheorem{lemma}[theorem]{Lemma}
\newtheorem{prop}[theorem]{Proposition}

\theoremstyle{definition}
\newtheorem{definition}[theorem]{Definition}
\theoremstyle{remark}
\newtheorem{example}[theorem]{Example}
\newtheorem{remark}[theorem]{Remark}

\title[Polynomial realizing $\rho_{E,n}(G_K)$]{Polynomials realizing images of Galois representations of an elliptic curve}
\author{Zoé Yvon}
\address{
Zoé Yvon, Institut de Mathématiques de Marseille (UMR 7373)
Site Sud, Campus de Luminy
Case 930
13288 MARSEILLE Cedex 9}
\email{zoenovy@free.fr}
\date{2022}

\newcommand\Aut{\operatorname{\mathrm{Aut}}}
\newcommand\Gal{\operatorname{\mathrm{Gal}}}

\newcommand\got\mathfrak
\newcommand{\N}{\mathbb{N}}
\newcommand{\Z}{\mathbb{Z}}
\newcommand{\Q}{\mathbb{Q}}
\newcommand{\F}{\mathbb{F}}

\newcommand{\id}{\operatorname{\mathrm{id}}}

\newcommand{\Norme}{\operatorname{\mathrm{N}}}

\usepackage{cancel}

\begin{document}

\begin{abstract}
The aim of the inverse Galois problem is to find extensions of a given field whose Galois group is isomorphic to a given group. In this article, we are interested in subgroups of $\mathrm{GL}_2(\Z/n\Z)$ where $n$ is an integer. We know that, in general, we can realize these groups as the Galois group of a given number field, using the torsion points on an elliptic curve. Specifically, a theorem of Reverter and Vila gives, for each prime $n$, a polynomial, depending on an elliptic curve, whose Galois group is $\mathrm{GL}_2(\Z/n\Z)$. In this article, we generalize this theorem in several directions, in particular for $n$ non necessarily prime.
We also determine a minimum for the valuations of the coefficients of the polynomials arising in our construction, depending only on $n$.
\end{abstract}
\maketitle

\section*{Introduction and Notations}

Let $K$ be a number field. Let $\overline{K}$ be an algebraic closure of $K$ and $G_K=\Gal(\overline{K}/K)$ its absolute Galois group. Let $E/K$ be an elliptic curve, with point at infinity $\mathcal{O}$. For a positive integer $n$, we consider the Galois representation $\rho_{E,n}$ corresponding to the action of $G_K$ on $E[n]$, the group of the $n$-torsion points of $E$, that is
\[\rho_{E,n}:G_K\to\Aut(E[n])\simeq\mathrm{GL}_2(\Z/n\Z),\]
following \cite[III.7]{AEC}.
The kernel of the Galois representation $\rho_{E,n}$ is the set of $\sigma\in G_K$ such that $\sigma$ acts trivially on the $n$-torsion points, in other words, $\sigma$ fixes the extension of $K$ generated by the coordinates of the $n$-torsion points, denoted by $K(E[n])$. So $\mathrm{ker}\rho_{E,n}=\Gal(\overline{K}/K(E[n]))$, and it is a normal subgroup of $G_K$. Therefore, the extension $K(E[n])/K$ is Galois and $\Gal(K(E[n])/K)\simeq\mathrm{Im}\rho_{E,n}$. This gives a solution for the inverse Galois problem over $K$ for the subgroups of $\mathrm{GL}_2(\Z/n\Z)$, which are images of such representations. We can call such a realisation \emph{geometric}, because it uses elliptic curves. More generally, we can talk about geometric realization when we use abelian varieties to realize groups as Galois groups. This method makes it possible to tabulate extensions of number fields with fixed Galois group and specified ramification. Indeed, the ramification of $K(E[n])/K$ is controlled by $n$ and the conductor of $E$.

Now we know that $\rho_{E,n}(G_K)$ is realizable as Galois group over $K$, we would like to have an explicit polynomial whose Galois group is $\rho_{E,n}(G_K)$, \emph{i.e.} whose roots generate $K(E[n])/K$. Reverter and Vila found such a polynomial for $n$ prime, and $\rho_{E,n}$ surjective. Their theorem is the starting point of this article:

\begin{theorem*}\label{Reverter-Vila}(\cite[Theorem 2.1]{Rvil}) Let $K$ be a field of characteristic $0$. Let $E$ be an elliptic curve over $K$, given by Weierstrass equation
$E:y^2=f(x)$, and let $\ell$ be an odd prime. Suppose that the Galois representation
\[\rho_{E,\ell}:G_K\to\Aut(E[\ell])\simeq\mathrm{GL}_2(\F_\ell)\]is surjective and let 
$P\in E[\ell]\setminus \{\mathcal{O}\}$. Then
\begin{enumerate}
\item The $\ell$-division polynomial is irreducible and its Galois group over $K$ is isomorphic to $\rho_{E,\ell}(G_K)/\{\pm\id\}$.
\item The characteristic polynomial $\chi_{E,\ell}$ of the multiplication by $x(P)+y(P)$ in $K(P):=K(x(P),y(P))$ is irreducible with Galois group isomorphic to $\rho_{E,\ell}(G_K)$.
\end{enumerate}
\end{theorem*}

The hypothesis "the representation $\rho_{E,\ell}$ is surjective" is not very restrictive since, in \cite{Serre}, Serre proved that, if $E$ does not have complex multiplication, then $\rho_{E,\ell}$ is surjective for almost all prime $\ell$. But it does not cover all of the cases. In the case where the ground field is $\Q$, these remaining cases are given by Zywina: in \cite{zywinasurj}, he gives an algorithm that exactly determines the exceptional primes for $E/\Q$, and in \cite{Zywinapassurj} he gives the possible images of non surjective representations for $E/\Q$. Also, in \cite{sutherland_2016}, Sutherland gave an algorithm to compute images mod $\ell$, for $\ell$ prime, for elliptic curves over $\Q$. For $E$ over a general number field, Larson and Vaintrob give, under the Generalized Riemann Hypothesis, many upper bounds for the exceptional primes, see \cite{10.1112/blms/bdt081}. 

In this article, we generalized Reverter and Vila's result for any equation for $E$, any image of $\rho_{E,n}$, with $n$ an integer not necessary prime, and we consider a broader choice of functions in $K(E)$ than $x+y$. Since $n$ is not necessarily prime, in Subsection~\ref{primitive division polynomials} we define the $n$-th primitive division polynomial $\widetilde{\psi}_n$, which is a factor of the $n$-th division polynomial corresponding to the points of exact order $n$. The main result of the first section of this article is the following theorem.
\begin{theorem*}(Theorems~\ref{galois psi},~\ref{pol réalisant im rho} and~\ref{transitive case})
Let $E$ be an elliptic curve over $K$ with Weierstrass equation $w_E(x,y)=0$. Let $u\in K(E)$ with degree $1$ in $x$ and $y$ such that $K(u,[-1]^*u)=K(x,y)$, and $n\geq3$.
\begin{enumerate}
    \item The $n$-th primitive division polynomial $\widetilde{\psi}_n$ has Galois group isomorphic to $\rho_{E,n}(G_K)/\{\pm\id\}$. If the action of $G_K$ on the points of order $n$ is transitive, then $\widetilde{\psi}_n$ is irreducible.
    \item The characteristic polynomial $\chi_{u,n}$ of multiplication by $u$ in the ring $K[X,Y]/(w_E,\widetilde{\psi}_n)$ has Galois group isomorphic to $\rho_{E,n}(G_K)$. Moreover, $\chi_{u,n}$ is irreducible if and only if $G_K$ acts transitively on the points of order $n$.
\end{enumerate}
\end{theorem*}

The first point of the theorem is well known and pretty immediate, the novelty lies in the second point.

The degree of $\chi_{u,n}$ is $2\deg\widetilde{\psi}_n$ which is less than or equal to $n^2-1$. So this theorem gives a way to construct polynomials of high degree with known Galois groups, which are subgroups of $\mathrm{GL}_2(\Z/n\Z)$, for an arbitrary integer $n$.

In Section~\ref{polynomial realizing im rho}, we prove this theorem step by step, together with a case study of the case $n=3$, in Subsection~\ref{case n=3}.
Section~\ref{valuation} focuses on the arithmetic properties of the coefficients of $\chi_{u,n}$. First, we give a optimal lower bound for their valuations in Proposition~\ref{minimum valuation}, then discuss the computation of $\chi_{u,n}$, see Remark~\ref{two ways to calculate}. Section~\ref{examples} contains some examples, one of them, Example~\ref{ex serre curves}, based on a family of Serre curves given by \cite{Dserrecurves}.

Most of results are valid in more general fields, but we restrict our interests to number fields.

In this article we will use the following notations: for a polynomial $f$ in $K[X]$, we denote by $K(f)\subset\overline{K}$ the splitting field of $f$ over $K$, that is the extension of $K$ generated by the roots of $f$, and by $\Gal(f)$ the Galois group of $K(f)/K$. For $u\in K(E)$ and $A$ a subgroup of $E(\overline{K})$, we denote by $u(A)$ the set of the $u(P)$ such that $P\in A$.

\section*{Acknowledgement}

I thank my supervisors Samuele Anni and David Kohel.

\section{Polynomial realizing $\rho_{E,n}(G_K)$}\label{polynomial realizing im rho}

Let $E/K$ be an elliptic curve. For a positive integer $n$, we define \[E_n:=\{P\in E(\overline{K}) \text{ of order }n\}\subset E[n].\]

\subsection{Primitive division polynomials}\label{primitive division polynomials}

Let $(\psi_n)$ be the family of the $n$-division polynomials of $E$, for $n\in\N$; for references see for example \cite[Exercise 3.7]{AEC} and \cite{McKee1994ComputingDP}.
Studying the points of order $n$, where $n$ is not necessary prime, naturally leads to the following definition:

\begin{definition}Let $E$ be an elliptic curve and $(\psi_n)$ the family of its division polynomials. We define the \emph{primitive division polynomial} $(\widetilde{\psi}_n)$ recursively by\begin{equation*}\psi_n=\underset{m\mid n}{\prod}\widetilde{\psi}_m.\end{equation*}
\end{definition}

\begin{remark}\label{pol primitif et pol div}
Clearly $\psi_1=\widetilde{\psi}_1=1$ and, for $p$ prime and $k\geq1$, we have $\widetilde{\psi}_{p^k}=\frac{\psi_{p^k}}{\psi_{p^{k-1}}}$. In particular, $\widetilde{\psi}_p=\psi_p$ for $p$ prime.
\end{remark}

\begin{prop}\label{action sur ordre n}If $n\neq2$, the polynomial $\widetilde{\psi}_n$ is in $K[X]$ and its roots are the elements of $x(E_n)$. Moreover, if $n$ is a power of a prime $p$ then the leading coefficient of $\widetilde{\psi}_n$ is $p$. Otherwise, the polynomial $\widetilde{\psi}_n$ is monic.
\end{prop}

\begin{proof}
For $n$ odd, respectively even, the polynomial $\psi_n$, respectively $\frac{\psi_n}{\psi_2}$, is univariate with roots the elements of $x(E[n])$, respectively $x(E[n])\setminus x(E[2])$; see \cite[exercice 3.7.(a) and (d)]{AEC}. Then, by induction, for $n\geq3$, the polynomial $\widetilde{\psi}_n$  is univariate with roots the elements of $x(E_n)$.
Now, the absolute Galois group $G_K$ acts on the $n$-torsion points, therefore on the points of order $n$. Indeed, let $P$ be a point of order $n$ and $\sigma$ in $G_K$. Suppose that $\sigma(P)$ has order $m<n$. Then, the point $P=\sigma^{-1}(\sigma(P))$ belongs to $E[m]$, which is a contradiction. Therefore, the factorisation \[\psi_n=\underset{m\mid n}{\prod}\widetilde{\psi}_m\]is defined over $K$. For a polynomial $f$, let $c(f)$ be its leading coefficient. We have \[c(\psi_n^2)=\underset{m\mid n}{\prod}c(\widetilde{\psi}_m^2).\]
So, using that $c(\psi_n^2)=n^2$ for all $n$ and Remark~\ref{pol primitif et pol div}, we conclude, by induction that $c(\widetilde{\psi}_{p^k})=p$ for $p$ prime and $k$ a positive integer, and so that $c(\widetilde{\psi}_n)=1$ if $n$ is composite.
\end{proof}

\begin{remark}
Since we know the degree of $\psi_m^2$ for all $m$, we can compute the degree of $\widetilde{\psi}_n$, as a polynomial in $x$, by induction. For example, if $n$ is the product of two distinct primes $p$ and $q$, then $\widetilde{\psi}_n$ has degree $(n^2-p^2-q^2+1)/2$.
\end{remark}

The polynomials $\psi_n$ and $\psi_2\psi_n$ are in $K[x]$ when $n$ is odd and even respectively, and $\psi_2$ is not. In particular, $K(\psi_n)$ is well-defined when $n$ is odd, but not when $n$ is even.

\begin{definition}
For $n$ an even integer, we define $K(\psi_n):=K(\psi_2\psi_n)$.
\end{definition}

\begin{lemma}For $n$ a positive integer, we have \[K(\psi_n)=K(x(E[n]))=K(x(E_n))=K(\tilde{\psi_n}).\]
\end{lemma}

\begin{proof}
The last equality is given by Proposition~\ref{action sur ordre n}.

By Proposition~\ref{action sur ordre n} and the factorization of $\psi_n$, the roots of $\psi_n$, for $n$ is odd, are the elements of $x(E[n])$. The same is true for $\psi_2\psi_n$ for $n$ even, noting that $x(E[2])$ are the roots of $\psi_2^2$, see \cite[exercice 3.7.(d)]{AEC}. Then we have the first equality.

Finally, for the second equality, since $E_n\subset E[n]$, then $\widetilde{\psi}_n$ divides $\psi_n$, and we obviously have $K(\widetilde{\psi}_n)\subset K(\psi_n)$.
For the reverse inclusion, let $x(P)\in K(x(E[n]))$ with $P$ of order $m\mid n$. Then $n=km$ for some $k$, and $P=kQ$ for some point $Q$ of order $n$. From \cite[exercise 3.7.(d)]{AEC}, \[x(P)=x(kQ)=\frac{\phi_k(x(Q))}{\psi_k^2(x(Q))},\]where $\phi_k\in K[X]$. So $x(P)\in K(x(E_n))=K(\widetilde{\psi}_n)$.
\end{proof}

In fact, we do not need these formulas to prove that $K(x(E[n]))=K(x(E_n))$: see Lemma~\ref{égalité extension torsion et ordre} for a proof using only Galois theory.

\subsection{Galois group of $\psi_n$}\label{groupe de galois de psi_n}

For a positive integer $n$, let $\pi_n$ be the canonical projection \[\pi_n:\mathrm{GL}_2(\Z/n\Z)\to \mathrm{GL}_2(\Z/n\Z)/\{\pm\id\},\]and define $\overline{\rho_{E,n}}=\pi_n\circ\rho_{E,n}$. If $\rho_{E,n}(G_K)$ does not contain $-\id$, it is canonically isomorphic to $\overline{\rho_{E,n}}(G_K)$.

\begin{lemma}\label{opposé y coordonnée}For $n\geq2$, we have $\ker\overline{\rho_{E,n}}=\Gal(\overline{K}/K(x(E_n)))$.
\end{lemma}

\begin{proof}
If $\sigma\in G_K$ satisfies $\overline{\rho_{E,n}}(\sigma)=\id$, then $\sigma(x(P))=x(P)$ for all $P\in E_n$.
Now, let $\sigma\in \Gal(\overline{K}/K(x(E_n)))$. For $P\in E_n$, we have $\sigma(P)\in\{\pm P\}$. Suppose that there are two points $P,Q\in E_n$ such that $\sigma(P)=-P$ and $\sigma(Q)=Q$. Then, on the one hand,
\[\sigma(P+Q)=-P+Q\]and, on the other hand,
\[\sigma(P+Q)\in\{\pm (P+Q)\}.\] So either $P$ or $Q$ has order $2$, hence $n=2$, and in this case $P=-P$ for all $P\in E_n$.
Consequently, either $\sigma(P)=P$ for all $P\in E_n$, or $\sigma(P)=-P$ for all $P\in E_n$.
\end{proof}

\begin{lemma}\label{égalité extension torsion et ordre}For $n\geq2$, we have $K(x(E_n))=K(x(E[n]))$. \end{lemma}

\begin{proof}
The inclusion $K(x(E_n))\subset K(x(E[n]))$ is obvious. Then we have \[\Gal(\overline{K}/K(x(E[n]))<\Gal(\overline{K}/K(x(E_n)).\]Now, take $\sigma\in \Gal(\overline{K}/K(x(E_n)))$. By Lemma \ref{opposé y coordonnée}, we have $\overline{\rho_{E,n}}(\sigma)=\id$, hence $\sigma$ fixes $K(x(E[n]))$. 
\end{proof}

\begin{theorem}\label{galois psi}Let $n\geq2$.
\begin{enumerate}
\item The Galois group of $K(\widetilde{\psi}_n)$ over $K$ is isomorphic to $\overline{\rho_{E,n}}(G_K)$.
\item If $\rho_{E,n}(G_K)$ contains $-\id\neq\id$, then the extension $K(E[n])/K(x(E[n]))$ has degree $2$. Otherwise, $K(E[n])=K(x(E[n]))$ and the Galois group of $K(\widetilde{\psi}_n)$ is isomorphic to $\rho_{E,n}(G_K)$.
\end{enumerate}
\end{theorem}

\begin{remark}
\cite[Theorem 1.1]{Rvil} also gives a sufficient but not necessary criterion for the Galois group of $\widetilde{\psi}_n$ to be $\rho_{E,n}(G_K)$. This criterion implies in particular that $-\id$ is in the image, but this is not an equivalence.
\end{remark}

\begin{proof}
Using Lemma \ref{opposé y coordonnée}, we obtain \[\overline{\rho_{E,n}}(G_K)\simeq G_K/\ker\overline{\rho_{E,n}}\simeq\Gal(K(x(E_n))/K)=\Gal(K(\widetilde{\psi}_n)/K).\]
The second point of the proposition is also an immediate consequence of Lemma~\ref{opposé y coordonnée}.
\end{proof}

\subsection{Consequence in the case $n=3$}\label{case n=3}

Let $E/K$ be an elliptic curve and let $G:=\rho_{E,3}(G_K)$. We denote by $\zeta_3$ a primitive third root of unity.

We recall the following result. Let $n$ be a positive integer and $\zeta_n$ be a primitive $n$-th root of unity. By the Weil pairing, we have $K(\zeta_n)\subset K(E[n])$. For each $\sigma\in \Gal(K(E[n])/K)$, there exists an $\alpha(\sigma)\in (\Z/n\Z)^*$ such that $\sigma(\zeta_n)=\zeta_n^{\alpha(\sigma)}$. Thanks to the Weil pairing again, $\alpha(\sigma)$ satisfies $\det\circ\rho_{E,n}(\sigma)=\alpha(\sigma).$
Then the image of $\det\circ\rho_{E,n}$ is equal to the image of the natural embedding of $\Gal(K(\zeta_n)/K)$ in $(\Z/n\Z)^*$, which is surjective if and only if $K$ does not contain any $n$-th roots of unity.
In the case $n=3$, it gives that $\det\circ\rho_{E,3}$ has image $\{\pm1\}$ if and only if $K$ does not contain $\zeta_3$.

\begin{figure}[ht!]\caption{Subgroup lattice of $\mathrm{GL}_2(\F_3)$\label{graph 3}}\[\scalebox{0.60}{\begin{tikzpicture}[
roundnode/.style={circle, draw=green!60, fill=green!5, very thick, minimum size=7mm},
squarednode/.style={rectangle, draw=red!60, fill=red!5, very thick, minimum size=5mm},
]
%Nodes
\node      (GL) {$\mathrm{GL}_2(F_3)$};
\node   (N)       [below=of GL] {$\underset{\displaystyle \widetilde{D}_{16}}{\left\langle\begin{pmatrix}1&-1\\1&1\end{pmatrix},\begin{pmatrix}1&0\\0&-1\end{pmatrix}\right\rangle}$};
\node (SL) [left=of N]{$\mathrm{SL}_2(\F_3)$};
\node      (B)       [right=of N] {$\underset{\displaystyle D_{12}}{\begin{pmatrix}*&*\\0&*\end{pmatrix}}$};
\node        (C)       [below=of SL] {$\underset{\displaystyle C_8}{\left\langle\begin{pmatrix}  1&-1\\1&1\end{pmatrix} \right\rangle}$};
\node (Q) [left=of C] {$\underset{\displaystyle Q_8}{\left\langle \begin{pmatrix}0&1\\-1&0\end{pmatrix},\begin{pmatrix} 1&1\\1&-1\end{pmatrix}\right\rangle}$};
\node (D) [right=of C] {$\underset{\displaystyle D_8}{\begin{pmatrix}*&0\\0&*\end{pmatrix}\cup\begin{pmatrix}0&*\\ *&0\end{pmatrix}}$};
\node (S31) [below=of B] {$\underset{\displaystyle 1S_3}{\begin{pmatrix}1&*\\0&*\end{pmatrix}}$};
\node (S32) [right=of S31] {$\underset{\displaystyle 2S_3}{\begin{pmatrix}*&*\\0&1\end{pmatrix}}$};
\node (C6) [right=of S32] {$\underset{\displaystyle C_6}{\left\langle\begin{pmatrix}
-1&-1\\0&-1
\end{pmatrix}\right\rangle}$};
\node (C4) [below=of C] {$\underset{\displaystyle C_4}{\left\langle\begin{pmatrix}
0&1\\-1&0
\end{pmatrix}\right\rangle}$};
\node (D4) [below=of D] {$\underset{\displaystyle V_4}{\begin{pmatrix}*&0\\0&*\end{pmatrix}}$};
\node (-id) [below=of D4] {$\underset{\displaystyle 1C_2}{\langle -\id\rangle}$};
\node (C2) [right=of -id] {$\underset{\displaystyle 2C_2}{\begin{pmatrix}1&0\\0&*\end{pmatrix}}$};
\node (C3) [right=of C2] {$\underset{\displaystyle C_3}{\left\langle\begin{pmatrix}
1&1\\0&1
\end{pmatrix}\right\rangle}$};
\node (id) [below=of C2] {$\{\id\}$};

%Lines
\draw[-] (S32.south) -- (C3.north);
\draw[-] (S31.south) -- (C3.north west);
\draw[-] (GL.south west) -- (SL.north);
\draw[-] (GL.south) -- (N.north);
\draw[-] (GL.south east) -- (B.north);
\draw[-] (SL.south west) -- (Q.north);
\draw[-] (N.west) -- (Q.north east);
\draw[-] (N.south west) -- (C.north);
\draw[-] (N.south) -- (D.north);
\draw[-] (B.south west) -- (D4.north east);
\draw[-] (B.south) -- (S31.north);
\draw[-] (B.south east) -- (S32.north);
\draw[-] (B.east) -- (C6.north west);
\draw[-] (Q.south) -- (C4.north west);
\draw[-] (C.south) -- (C4.north);
\draw[-] (D.south west) -- (C4.north east);
\draw[-] (D.south) -- (D4.north);
\draw[-] (C4.south) -- (-id.north west);
\draw[-] (D4.south) -- (-id.north);
\draw[-] (D4.south east) -- (C2.north west);
\draw[-] (S31.south west) -- (C2.north);
\draw[-] (S32.south west) -- (C2.north east);
\draw[-] (C6.south west) -- (C3.north east);
\draw[-] (-id.south) -- (id.north west);
\draw[-] (C2.south) -- (id.north);
\draw[-] (C3.south) -- (id.north east);
\draw[-] (C6.south) .. controls +(down:70mm) and +(down:35mm) .. (-id.south west);
\draw[-] (SL.north west) .. controls +(up:45mm) and +(up:40mm) .. (C6.north);
\end{tikzpicture}}.\]\end{figure}
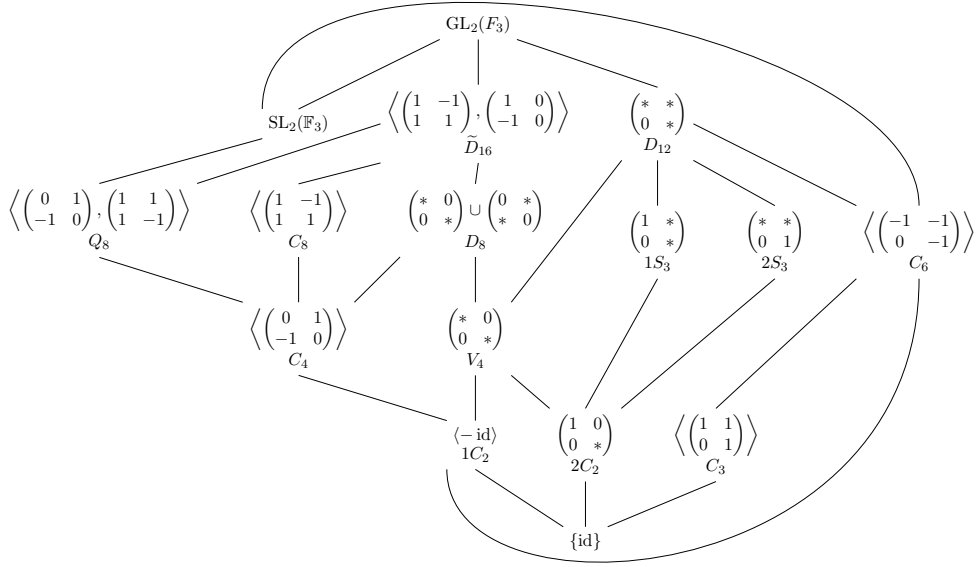

Since $G$ is a subgroup of $\mathrm{GL}_2(\F_3)$, we give a classification of all subgroups of $\mathrm{GL}_2(\F_3)$, up to conjugation, in Figure~\ref{graph 3}. Here, $C_n$ is the cyclic group of order $n$, $S_n$ is the symmetric group of degree~$n$, $D_{2n}$ is the dihedral group of order $2n$, $V_4$ is the Klein group, $Q_8$ is the quaternion group and $\widetilde{D}_{16}$ is the quasi-dihedral group of order $16$. The groups $1S_3$ and $2S_3$ are both isomorphic to $S_3$ but are not conjugate, and similarly for $1C_2$ and $2C_2$, both isomorphic to $C_2$.
The stars $*$ in the matrices means that all choices of elements in $\F_3$ are possible, provided the matrix is invertible. The graph shows directly
\begin{itemize}
    \item which subgroups contain $-\id$, so, conversely, which images are isomorphic to the Galois group of $\psi_3$, 
    \item and which subgroups are in $\mathrm{SL}_2(\F_3)$ or not, which only depends on whether $K$ contains $\Q(\zeta_3)$.
\end{itemize}

\begin{remark}
Following the classical description of subgroups of $\mathrm{GL}_2(\Z/n\Z)$ given by Serre in \cite[Section 2]{Serre} or Sutherland in \cite{sutherland_2016}, $\begin{pmatrix}*&0\\0&*\end{pmatrix}$ is $C_s(3)$, the split Cartan subgroup; $\begin{pmatrix}*&0\\0&*\end{pmatrix}\cup\begin{pmatrix}0&*\\ *&0\end{pmatrix}$ is $N_s(3)$, the normalizer of $C_s(3)$; $\left\langle\begin{pmatrix}  1&-1\\1&1\end{pmatrix} \right\rangle$ is $C_{ns}(3)$, the non split Cartan subgroup; $\left\langle\begin{pmatrix}1&-1\\1&1\end{pmatrix},\begin{pmatrix}1&0\\0&-1\end{pmatrix}\right\rangle$ is $N_{ns}(3)$, the normalizer of $C_{ns}(3)$;
and $\begin{pmatrix}*&*\\0&*\end{pmatrix}$ is $B(3)$, the full Borel subgroup.
\end{remark}

The next two propositions give all possibilities for $G$, the image of the mod~$3$ representation, depending on the factorization of $\psi_3$.

\begin{theorem}\label{im mod 3 1}
If $K\cap \Q(\zeta_3)=\Q$, then, using the classification above,
\begin{enumerate}
    \item If $\psi_3$ is irreducible, then\begin{enumerate}\item $G=\mathrm{GL}_2(\F_3)$ if $K(\psi_3)/K$ is generated by exactly three roots of $\psi_3$; 
    \item $G\simeq \widetilde{D}_{16}$ if $K(\psi_3)/K$ is generated by exactly two roots of $\psi_3$; or  
    \item $G\simeq C_8$ if $K(\psi_3)/K$ is cyclic.
\end{enumerate}
    \item If $\psi_3$ can be factorized in two irreducible polynomials of degree $2$, then $G\simeq D_8$.
    \item If $\psi_3$ has a unique root over $K$ then \begin{enumerate}
        \item $G\simeq D_{12}$  if $K(E[3])/K(x(E_3))$ has degree $2$; or
\item $G\simeq 1S_3$ or $2S_3$ if $K(E[3])=K(x(E_3))$.
    \end{enumerate}
\item If $\psi_3$ has exactly two roots over $K$, then
\begin{enumerate}
\item  $G\simeq V_4$ if $K(E[3])/K(x(E[3]))$ has degree $2$; or 
\item  $G\simeq 2C_2$ if $K(E[3])=K(x(E[3]))$.
\end{enumerate}
\end{enumerate}
\end{theorem}

\begin{theorem}\label{im mod 3 2}
If $K\cap \Q(\zeta_3)\neq\Q$, then, using the classification above,
\begin{enumerate}
    \item If $\psi_3$ is irreducible, then\begin{enumerate}
        \item $G=\mathrm{SL}_2(\F_3)$; or
        \item  $G\simeq Q_8$ if $K(\psi_3)/K$ is cyclic.
    \end{enumerate}
\item If $\psi_3$ factors into two irreducible polynomials of degree $2$ then $G\simeq C_4$.
    \item If $\psi_3$ has a unique root in $K$, then 
\begin{enumerate}
\item $G\simeq C_6$ if $K(E[3])/K(x(E[3]))$ has degree $2$; or
\item $G\simeq C_3$ if $K(E[3])=K(x(E[3]))$.
\end{enumerate}
\item If $\psi_3$ splits completely over $K$ then
\begin{enumerate}
\item $G\simeq C_2$ if $K(E[3])/K(x(E[3]))$ has degree $2$; or
\item $G\simeq C_1$ if $K(E[3])=K(x(E[3]))$.
\end{enumerate}
\end{enumerate}
\end{theorem}

\begin{proof}(of Theorems~\ref{im mod 3 1} and~\ref{im mod 3 2}). We will consider the Galois action on the $x$-coordinates on the $3$-torsion points, that are the roots of $\psi_3$, as permutations on these roots. Each $\sigma\in G_K$ corresponds to an element of $S_4$, seen as permutation representation, and is mapped to $\overline{\rho_{E,3}}(\sigma)$ in $\overline{G}:=G/\{\pm\id\}$. 
So $\overline{G}\simeq\Gal(\psi_3)$ corresponds to a subgroup of $S_4\simeq \mathrm{Sym}(\{\text{roots of }\psi_3\})$. Then, $G\subset B(3)$ is equivalent to $\psi_3$ having a root in $K$.
The action of the groups not in $B(3)$ or $N_s(3)$ is transitive on $x(E_3)$, so $\psi_3$ is irreducible. Consequently, $G\subset N_s(3)$ is equivalent to $\psi_3$ factoring into two degree $2$ polynomials over $K$.
Moreover, when $G$ does not contain $\{-\id\}$,
then $\Gal(\psi_3)\simeq G$, by Thereom~\ref{galois psi}. Otherwise, $\Gal(\psi_3)$ has index $2$ in $G$.

If $K\cap\Q(\zeta_3)=\Q$, we consider only the groups in the above graph with surjective determinant, \emph{i.e.}~not contained in $\mathrm{SL}_2(\F_3)$. We just observed that, if $\psi_3$ splits over $K$, then $G\subset\{\pm\id\}\subset\mathrm{SL}_2(\F_3)$. Hence, this cannot happen if $K\cap\Q(\zeta_3)=\Q$.

If $\Q(\zeta_3)\subset K$, we only consider the groups in the above graph with non surjective determinant, \emph{i.e.}~contained in $\mathrm{SL}_2(\F_3)$. If $\psi_3$ has at least two roots in $K$, then there is a basis of $E[3]$ such that $G\subset C_s(3)$. Since $\Q(\zeta_3)\subset K$, the group $G$ has determinant $1$, so $G\subset\{\pm\id\}$, hence $\psi_3$ splits over $K$.

The different items of each proposition are immediate deductions from these observations.
\end{proof}

\begin{corollaire}
If $\psi_3$ is irreducible or if $\psi_3$ factors into two irreducible polynomials of degree $2$, then $-\id$ belongs to the image.
\end{corollaire}

\subsection{Polynomials generating the image of $\rho_{E,n}$}

If $m$ divides $n$, the coordinates of points of order $m$ on an elliptic curve $E/K$ are obtained by adding $n/m$ times a point of order $n$. The addition of points is given by rational functions over $K$, then we have \[K(E[n])=K(E_n)\]for all positive integers $n$.

\subsubsection{General settings}

Let $E/K$ be an elliptic curve defined by a Weierstrass equation $w_E(x,y)=0$, so we have $K(E)=K(x,y)$.
Let $u$ be a polynomial of degree $1$ in $x$ and $y$ such that \begin{equation}\label{condition for u}K[u,u^*]=K[x,y],\tag{$*$}\end{equation} with $u^*:=[-1]^*u$. In particular, we have $u\neq u^*$. This will be our recurring hypothesis on $u$.\bigskip

Let $n$ be a positive integer. Let \[A=\left( K[X,Y]/(\widetilde{\psi}_n(X,Y),w_E(X,Y)\right)\simeq\left\{\begin{array}{ll} \frac{K[X]}{(\widetilde{\psi}_n)}\oplus \frac{K[X]}{(\widetilde{\psi}_n)}y &\mbox{if } n\neq2\\\frac{K[X]}{(\psi_2^2(X))}&\mbox{if } n=2.\end{array}\right.\]
The dimension of $A$ over $K$ is finite, and equal to $2\deg(\widetilde{\psi}_n)$ for $n\neq2$ and to $3$ for $n=2$, that is, in all cases, equal to cardinality $\vert E_n\vert$. 
We denote by $\chi_{u,n}$ the characteristic polynomial of the multiplication by the polynomial $u$ on the ring $A$. It is a monic polynomial of degree $\vert E_n\vert$ with coefficients in $K$. 
We denote by $g_1,\dots, g_s\in K[X]$ the irreducible and monic factors of $\widetilde{\psi}_n$, if $n\neq2$, or of $\psi_2^2$ if $n=2$.
The polynomials $g_i$ are coprime because the roots of $\psi_2^2$ and of $\widetilde{\psi}_n$ for $n\neq2$ all have multiplicity one.
We set $A_i=K[X,Y]/(g_i,w_E)$ and denote by $\chi_i$ the characteristic polynomial of the class of $u$ in $A_i$.
From the Chinese remainder theorem, we have \[A\simeq A_1\times \dots\times A_s.\]

\begin{lemma}\label{racines de chii}
The roots of $\chi_i$ are $u(P)$ where $P\in E_n$ and $g_i(x(P))=0$.
\end{lemma}

\begin{proof}
Let $S=\{P\in E_n\mid g_i(x(P))=0\}$.
We know that the roots of $\chi_i$ are in $u(S)$.
We have to prove the reverse inclusion. Let $P\in S$. If $n=2$ then $u(P)\in K(x(P))$. Therefore, since the elements of $x(S)$ are conjugate, the elements of $u(S)$ are too. Now, suppose $n\geq3$.
The polynomial $w_E$ has degree $2$ over the field $K[X]/(g_i)$. If $w_E$ is irreducible over $K[X]/(g_i)$, take $P'\in S$ such that $u(P')$ is a root of $\chi_i$. The irreducibility of $g_i$ implies that there exists $\sigma\in G_K$ sending $x(P')$ to $x(P)$, and so $u(P')$ to $u(P)$ or $u(-P)$. The irreducibility of $w_E$ over $K(x(P))$ implies that there exists $\tau\in G_K$ which fixes $K(x(P))$ and sends $y(P)$ to $y(-P)$. Therefore $u(P)$ is conjugates to $u(P')$, so it is a root of $\chi_i$.

If $w_E$ is not irreducible, then $w_E=(Y-\alpha)(Y-\beta)$ with $\alpha, \beta\in K[X]/g_i(X)$. So $A_i$ is the product of two extensions of $K$, that is
\[A_i\simeq (K[X]/g_i(X))[Y]/(Y-\alpha)\times(K[X]/g_i(X))[Y]/(Y-\beta).\]
And the result follows since, if $P\in S$, then $y(P)$ is $\alpha(x(P))$ or $\beta(x(P))$ so $u(P)$ is a root of $\chi_i$.
\end{proof}

\begin{prop}\label{roots of chi} We have
\[\chi_{u,n}(T)=\underset{P\in E_n}{\prod}(T-u(P)).\]
\end{prop}

\begin{proof}
From the isomorphism of $K$-vector spaces:\[A\simeq A_1\times\dots\times A_s,\]
the matrix of the multiplication by $u$ in $A$ is conjugated to a block matrix in the $K$-vector space $A_1\times\dots\times A_s$. The characteristic polynomials of these two matrices are the same, and each $\chi_i$ is the characteristic polynomial of the block corresponding to $A_i$. Therefore $\chi_{u,n}$ is the product of characteristic polynomials $\chi_i$ for $i=1,\dots,s$. The result follows from Lemma~\ref{racines de chii}.
\end{proof}

\begin{corollaire}\label{pol réalisant im rho}
We have $K(\chi_{u,n})\simeq K(E[n])$ and $\Gal(\chi_{u,n})\simeq\rho_{E,n}(G_K)$.
\end{corollaire}

\begin{proof}By Proposition~\ref{roots of chi} and assumption on $u$, we have
\begin{align*}K(\chi_{u,n})=K(u(E_n))=K(u(E_n),u^*(E_n))&=K(x(E_n),y(E_n))\\&=K(E_n)=K(E[n])\end{align*}
The group isomorphism follows.
\end{proof}

\subsubsection{Specialization}\label{Specialization}

If $A_i$ is a field, then \begin{itemize}
\item $A_i\simeq K(P)$ for $P\in E_n$ such that $g_i(x(P))=0$. 
\item from Lemma~\ref{racines de chii}, the $u(P)$ with $P\in E_n$ and $g_i(x(P))=0$ are all conjugate.
\end{itemize}
Therefore, the characteristic polynomial of $u(P)$ in $K(P)$ does not depend on the choice of $P\in E_n$ such that $g_i(x(P))=0$ and it is equal to $\chi_i$. If $A$ is a field, then $A\simeq K(P)$ for all $P\in E_n$ and the characteristic polynomial of $u(P)$ in $K(P)$ is $\chi_{u,n}$.
If $-\id$ belongs to $\rho_{E,n}(G_K)$ then $A_i$ is a field for all $i$. If the action of $G_K$ on $E[n]$ is transitive, then $A$ is a field.

\begin{prop}
Suppose that $\rho_{E,n}(G_K)$ contains $-\id$.
Then the splitting field of $\chi_i$ over $K$ is the extension generated by the roots of $g_i$ and their $y$-coordinates. The compositum of the splitting fields of the $\chi_i$ is $K(E[n])$, and its Galois group is isomorphic to $\rho_{E,n}(G_K)$.
\end{prop}

\begin{proof}Let $S=\{P\in E_n\mid g_i(x(P))=0\}$. We have \[K(\chi_i)=K(u(S),u^*(S))=K(x(S),y(S)).\]
Since $K(a,b)=K(a)K(b)$, the compositum of the extensions generated by the $\chi_i$ is the extension of $K$ generated by the roots of $\widetilde{\psi}_n$ if $n\geq3$ or of $\psi_2^2$, if $n=2$, and the corresponding $y$-coordinates, in other words the coordinates of the points of order $n$.
\end{proof}

\subsubsection{Condition on $u$}\label{condition on u}

The function $u$ is linear in $x$ and $y$ so it has the form $ay+bx+c$ with $a,b,c\in K$. What are the conditions on $a,b,c$ in terms of $E$ to have $K(u,u^*)=K(x,y)$?

Let $(a_i)$ be the coefficients of $E$. If $a\neq0$ and $2b-a_1a\neq0$ then \[x=\frac{u+u^*+aa_3-2c}{2b-a_1a}\]and \[y=\frac{(b-a_1a)u-bu^*-bab-3+a_1ac)}{(2b-a_1a)a}.\]
So we have $K(u,u^*)=K(x,y)$. We know that the condition $a\neq0$ is necessary. But, what about the condition $2b-a_1a\neq0$? Can we also take $(a,b)$ such that $2b-a_1a=0$? The following example show that this cannot be always the case.
\begin{example}Let $E$ be an elliptic curve of $j$-invariant $0$, with Weierstrass equation $y^2=x^3-B$, with $B\in K\setminus(K^*)^3$. Let $u=y$. In particular, with the previous notations, $a\neq 0$, $b=0$, $a_1=0$ and, then, $2b-a_1a=0$. Then $K(u,u^*)=K(y)$. The elliptic curve $E$ admits an automorphism given by $(x,y)\mapsto (\zeta_3 x,y)$. The required condition on $B$ makes the polynomial $x^3-B$ irreducible over $K$, then $\zeta_3x$ is conjugate to $x$ over $K(y)$, and $K(x,y)/K(y)$ has degree $3$, and not $1$. Hence $K(u,u^*)\neq K(x,y)$.
\end{example}

The required condition for $u$ is $K(u,u^*)=K(x,y)$. Actually, we can assume a weaker condition to obtain $\Gal(\chi_{u,n})=\rho_{E,n}(G_K)$. For example, if $-\id$ is not in $ \rho_{E,n}(G_K)$, then we are not forced to take $a\neq0$. Indeed, we can take $u=x$ and then $\chi_{u,n}$ is a scalar multiple of $\psi_n^2$, and its Galois group is $\rho_{E,n}(G_K)$.

\subsection{Criterion to have $-\id$ in the image}

The case where the image contains $-\id$ or not are clearly distinguished. So we can ask under which conditions this happens. Sutherland gives a criterion, linked with quadratic twists of the elliptic curve, in \cite[Lemma 5.24]{sutherland_2016} and additional conditions in \cite[Corollary 5.25]{sutherland_2016}.

Before giving an other criterion, here is a theorem of Serre which will be useful. For a prime ideal $\got{p}$ of the ring of integers $\mathcal{O}_K$ of a number field $K$, let $k_\got{p}$ be the residue field at $\got{p}$.

\begin{theorem}\label{Serre frob}(\cite[IV-5]{SerreMG})
Let $\got{p}$ be a prime ideal of $\mathcal{O}_K$ and $\ell$ be a prime such that $\ell$ does not divide $\mathrm{char}(k_\got{p})$. If $E$ has good reduction at $\got{p}$, then the Frobenius automorphism above $\got{p}$ is defined independently of the chosen prime ideal above $\got{p}$ and the characteristic polynomial of $\rho_{E,\ell}(\mathrm{Frob}_\got{p})$ is \[X^2-a_\got{p}(E)X+\Norme(\got{p})\pmod\ell\] where $a_\got{p}(E)=\Norme(\got{p})+1-\vert E_\got{p}(k_\got{p})\vert$ with $E_\got{p}$ the reduction of $E$ modulo $\got{p}$. 
\end{theorem}

We obtain the following proposition, with the same notations:

\begin{prop}Let $\ell$ be an odd prime. Let $E/K$ be an elliptic curve. The image $\rho_{E,\ell}(G_K)$ contains $-\id$ if and only if there exists a prime ideal $\got{p}$ of good reduction for $E$ such that $\Norme(\got{p})\equiv1\pmod\ell $ and $a_\got{p}(E)\equiv -2\pmod\ell$.
\end{prop}

\begin{proof}
Suppose that we have $\sigma\in G_K$ such that $\rho_{E,\ell}(\sigma)=-\id$. Then, from Cebotarev's density theorem, there exist infinitely many prime ideals $\got{p}$ of good reduction such that $\sigma=\mathrm{Frob}_\got{q}$ for a prime ideal $\got{q}$ of $\mathcal{O}_{\overline{K}}$ above $\got{p}$. Hence, the image of $\rho_{E,\ell}$ contains $-\id$ if and only if there exists a prime $\got{p}$ of good reduction such that $\rho_{E,\ell}(\mathrm{Frob}_\got{q})=-\id$. Since being conjugate to $-\id$ implies being equal to $-\id$, this equality is equivalent to the equality of the corresponding characteristic polynomials:\[X^2-a_\got{p}(E)X+\Norme(\got{p})\equiv X^2+2X+1\pmod\ell.\]The result follows.
\end{proof}

\begin{remark}
The proposition gives us a necessary and sufficient condition for $\rho_{E,\ell}(G_K)$ to contain $-\id$. In particular, it gives a criterion to know whether the image contains $-\id$, as required. In the other direction, if we know the image of the representation, it gives us a criterion on the cardinality of $E_\got{p}(k_\got{p})$ for a certain $\got{p}$.
\end{remark}

\begin{remark}
If $K$ is a number field, we also know that $\rho_{E,\ell}(G_K)$ is surjective for almost all primes $\ell$, and so almost always contains $-\id$. In the introduction, we gave some references about the non surjective cases.
\end{remark}

\begin{remark}
For a real $x$ and an integer $h$, Serre (\cite[Theorem 1.20]{Serre_Qlqaptdc}) gives an estimation of the number of primes $p\leq x$ such that $E/\Q$ has good reduction at $p$ and $a_p(E)=h$. This estimation is \[\mathcal{O}\left(\left(\frac{\mathrm{log}(x)}{\mathrm{log}\mathrm{log}(x)^2\mathrm{log}\mathrm{log}\mathrm{log}(x)}\right)^{1/4}\int_2^x \frac{dt}{\mathrm{log}(t)}\right),\] which does not depend on $h$. This theorem is valid also on a number field, by \cite[8.2, Remarques 2]{Serre_Qlqaptdc}.
\end{remark}

\begin{remark}Let $q:=\Norme(\got{p})$. Since we must have $q\equiv1\pmod\ell$, then $q$ is necessary bigger than $\ell$. If $\ell>2$, then $q\geq 2\ell+1$.
The Hasse-Weil bound \[-2\sqrt{q}+q+1\leq\vert E_\got{p}(k_\got{p})\vert\leq2\sqrt{q}+q+1.\]
imposes additional conditions on $\got{p}$.
\end{remark}

\subsection{Transitive case}

\begin{theorem}\label{transitive case}
Let $E$ be an elliptic curve over $K$ and let $n\geq2$. Suppose that the Galois action on $E_n$ is transitive.
Then
\begin{enumerate}
\item The polynomial $\widetilde{\psi}_n$, if $n\neq2$, and $\psi_2^2$, if $n=2$, is irreducible and its Galois group is isomorphic to  $\overline{\rho_{E,n}}(G_K)$.
\item The characteristic polynomial of the multiplication by $u(P)$ in $K(P)$ is irreducible and its Galois group is isomorphic to $\rho_{E,n}(G_K)$.
\end{enumerate}
\end{theorem}

\begin{proof}
Let $P,Q\in E_n$. Since the action is transitive, then $K[X,Y]/(w_E,\widetilde{\psi}_n)$ is a field, from Subsection~\ref{Specialization}. In particular, $\widetilde{\psi}_n$, if $n\neq2$, or $\psi_2^2$, if $n=2$, is irreducible and we know its Galois group from Theorem~\ref{galois psi}.

We know that $\Gal(\chi_{u,n})\simeq\rho_{E,n}(G_K)$, by Proposition~\ref{pol réalisant im rho} and Subsection~\ref{Specialization}. Therefore, we only have to show that $\chi_{u,n}$ is irreducible. In other words, it suffices to prove that its degree, $\vert E_n\vert$, is equal to the number of conjugate elements to $u(Q)$, for some $Q\in E_n$.
We have \[\{\sigma(u(Q)),\sigma\in G_K\}=\{u(P), P\in E_n\}.\] Let us show that the $u(P)$ are pairwise distinct. Suppose that $u(P)=u(P')$. Let $\sigma\in G_K$ be such that $\rho_{E,n}(\sigma)=-\id$. Then $\sigma(u(P))=\sigma(u(P'))$, in other words $u(-P)=u(-P')$ and \[u(P)\pm u(-P)=u(P')\pm u(-P').\] But, by assumption on $u$, we have $u(P)\neq u(-P)$, unless $P=-P$. Then, in all cases, $x(P)=x(P')$ and $y(P)=y(P')$, so $P=P'$. Therefore, the $u(P)$ are pairwise distinct, so their number is the cardinal of $E_n$.
\end{proof}

\begin{remark}\label{max image}
In particular, if $\rho_{E,n}$ is surjective, then $\widetilde{\psi}_n$ and $\chi_{u,n}$ are both irreducible and their Galois groups are respectively $\mathrm{GL}_2(\Z/n\Z)/\{\pm\id\}$ and $\mathrm{GL}_2(\Z/n\Z)$.
\end{remark}

\section{About the valuation of the coefficients of $\chi_{u,n}$}\label{valuation}

Let $E/K$ be an elliptic curve with Weierstrass equation $w_E(X,Y)=0$ where
\begin{equation}\label{equation E} w_E(X,Y)=Y^2+a_1XY+a_3Y-X^3-a_2X^2-a_4X-a_6.\tag{$\diamond$}\end{equation}Let $n$ be a positive integer. As in the first section, let $u\in K(E)$ be a function of degree $1$ in $x$ and $y$ such that $K(u,u^*)=K(x,y)$. We have seen in Subsection~\ref{condition on u} that $u=ay+bx+c$ for some $a,b,c\in K$ with $a\neq0$. Theorem~\ref{pol réalisant im rho} says that $\chi_{u,n}$ has Galois group $\rho_{E,n}(G_K)$. This second section gives a minimum for the valuation of the coefficients of $\chi_{u,n}$. As usual, the case $n=2$ has to be studied separately.

\subsection{Case $n=2$}

From Theorem~\ref{galois psi}, we have $\Gal(\psi_2^2)\simeq\rho_{E,2}(G_K)$. We can compute $\psi_2^2$:
\[\psi_2^2(x)=4x^3+(a_1^2+4a_2)x^2+(4a_4+2a_1a_3)x+a_3^2+4a_6.\]
But the polynomial $\psi_2^2$ is not normalized. We have \[\chi_{x,2}=\frac{1}{4}\psi_2^2=x^3+(\frac{a_1^2}{4}+a_2)x^2+(a_4+\frac{a_1a_3}{2})x+\frac{a_3^2}{4}+a_6\in \frac{1}{2^2}\Z[a_i][x].\]
If we have a short Weierstrass equation, that is $a_1=a_3=a_2=0$, then we obtain $E:y^2=\frac{1}{4}\psi_2^2$. Then a short Weierstrass equation for $E$ immediately gives a polynomial realizing $\rho_{E,2}(G_K)$. And the coefficients of $\psi_2^2$ have the smallest valuations possible at a prime ideal $\got{p}$ if we take a minimal Weierstrass equation for $E$ at $\got{p}$. The coefficients of $\psi_2^2$ have the smallest valuation possible at every prime ideal if we take a global minimal equation for $E$, that is possible, for example, if $K$ has class number $1$.

\subsection{Minimum of the valuations of the coefficients of $\chi_{u,n}$}

\begin{prop}\label{minimum valuation}Let $\ell$ be a prime and $k$ be an integer. Let $E/K$ be an elliptic curve given by~(\ref{equation E}).
Let $a,b,c\in K$, with $a\neq0$, and $R:=\Z[a_1,a_2,a_3,a_4,a_6]$. Then

\[\chi_{ay+bx+c,\ell^k}\in\frac{1}{\ell^3}R[a^{-1},b,c][X],\]
and, for $n$ a composite integer,
\[\chi_{ay+bx+c,n}\in R[a^{-1},b,c][X].\]

\end{prop}

\begin{proof}
We start by showing the result for $(b,c)=(0,0)$.
For $P\in E(\overline{K})$, let be the polynomial \[w_E(x(P),Y):=Y^2+(a_1x(P)+a_3)Y-(x(P)^3+a_2x(P)^2+a_4x(P)+a_6.\]
For $P\in E(\overline{K})$, the polynomial $w_E(x(P),Y)\in K(x(P))[Y]$ is monic, has degree $2$ and has roots $\pm y(P)$. If $P$ has order $n$, these are also roots of $\chi_{y,n}$. So $\underset{P\in E_n/\langle-\id\rangle}{\prod}w_E(x(P),Y)$ and $\chi_{y,n}(Y)$ are monic, both with degree equal to $2\deg\widetilde{\psi}_n$, and the same roots by Proposition~\ref{roots of chi}, so they are equal.
The resultant $\mathrm{Res}_X(\widetilde{\psi}_n,w_E)$, where $\widetilde{\psi}_n$ and $w_E$ are considered as polynomials in the first variable, belongs to $R[Y]$: it is a polynomial in $Y$ with coefficients in $R$. Let $r$ be the leading coefficient of $\widetilde{\psi}_n$. From \cite[A IV.75, Corollary 1]{Bourbaki2},  we have:\[\mathrm{Res}_X(\widetilde{\psi}_n,w_E)=r^{\deg_X(w_E)}\underset{\alpha\text{ roots of }\widetilde{\psi}_n}{\prod}w_E(\alpha,Y)=r^3\underset{P\in E_n/\langle{-\id}\rangle}{\prod}w_E(x(P),Y).\]
Therefore \[\chi_{y,n}\in\frac{1}{r^3}R[Y].\]
Now, consider the elliptic curve $E'/K$ which is the change of variables of $E$ given by 
$x=x'$ and $y=a^{-1}(y'-bx'-c)$. Then, we have a Weierstrass equation for $E'$, and $\chi_{y',n}$ has coefficients in $\frac{1}{r^3}R[a^{-1},b,c]$. So
\[\chi_{ay+bx+c,n}=\chi_{y',n}\in\frac{1}{r^3}R[a^{-1},b,c][Y].\]
By Proposition~\ref{action sur ordre n}, we obtain the desired result.
\end{proof}

\begin{remark}
The lower bound given in the proposition is a minimum, in the sense where there exists elliptic curves such that this bound is reached, as we will observe in Example~\ref{ex serre curves}.
\end{remark}

\begin{remark}\label{two ways to calculate}
The proof of the previous proposition points out the fact that we can compute $\chi_{u,n}$ in two ways. The first one is, by definition, to compute the characteristic polynomial of a matrix. The second one is, in view of the above, to compute the resultant of two polynomials. In the first case, we have to compute the determinant of a matrix of size $2\deg\widetilde{\psi}_n$, and in the second case the determinant of a matrix of size $\deg\widetilde{\psi}_n+3$. Since $\deg\psi_3=4$, the second matrix is always smaller, and the difference of size increases with $n$. Moreover, the matrix in the calculus of the resultant is easier to obtain, because we just have to put the coefficients of $\psi_n$ and $w_E$ and some zeros in the matrix. Whereas, to obtain the matrix of the multiplication by $u$, we have to choose a basis $(e_i)$ and write each $u*e_i$ in terms of the basis.

If $a_1\neq0$, then we can choose $u=y$ and $\chi_{y,n}$ is equal to
\[\mathrm{Res}_x(\widetilde{\psi}_n,w_E)=\det\begin{pmatrix}
r & 0 & 0 & 1 & 0 & \cdots & 0\\
\times & r & 0 & a_2 & \ddots & \ddots & \vdots\\
\vdots & \times & r & a_4-a_1y & \ddots & \ddots & 0 \\
\vdots & \vdots & \times & a_6-a_3y-y^2 & \ddots & \ddots & 1\\
\times & \vdots & \vdots & 0 & \ddots & \ddots & a_2\\
0 & \times & \vdots & \vdots & \ddots & \ddots & a_4-a_1y\\
0 & 0 & \times & 0 & \cdots & 0 & a_6-a_3y-y^2\end{pmatrix}\]
where $r$ and the crosses are elements of $K$, here corresponding respectively to the leading term of the polynomial $\widetilde{\psi}_n$ and its other coefficients.

If $a_1=0$, we can choose $u=x+y$ and $\chi_{x+y,n}$ is equal to \[\mathrm{Res}_{x'}(\widetilde{\psi}_n',w_{E'})=\det\begin{pmatrix}
r & 0 & 0 & 1 & 0 & \cdots & 0\\
\times & r & 0 & a_2-1 & \ddots & \ddots & \vdots\\
\vdots & \times & r & a_4-a_3-2y & \ddots & \ddots & 0 \\
\vdots & \vdots & \times & a_6-a_3y-y^2 & \ddots & \ddots & 1\\
\times & \vdots & \vdots & 0 & \ddots & \ddots & a_2-1\\
0 & \times & \vdots & \vdots & \ddots & \ddots & a_4-a_3-2y\\
0 & 0 & \times & 0 & \cdots & 0 & a_6-a_3y-y^2\end{pmatrix}\]where $E'$ is an elliptic curve obtained from $E$ setting by the change of variables $(x,y)\mapsto (x,x+y)$, and $\widetilde{\psi}_n'$ is its $n$-th primitive division polynomial. For comparison, here is the matrix of the multiplication by $y$,
\begin{equation*}\left(\begin{smallmatrix}
0&0&\cdots&\cdots&0&\times&a_6&0&\cdots&0&\times&\times&\times\\
0&0&\ddots&\ddots&\vdots&\times&a_4&a_6&\cdots&0&\times&\times&\times\\
0&0&\ddots&\ddots&\vdots&\vdots&a_2&a_4&\ddots&\vdots&\vdots&\vdots&\vdots\\
\vdots&0&\ddots&\ddots&\vdots&\vdots&1&a_2&\ddots&\vdots&\vdots&\vdots&\vdots\\
\vdots&\vdots&\ddots&\ddots&\vdots&\vdots&0&\ddots&\ddots&\vdots&\vdots&\vdots&\vdots\\
0&\vdots&\ddots&\ddots&1&\times&\vdots&\ddots&\ddots&1&\times&\times&\vdots\\
1&0&\ddots&\ddots&0&0&-a_3&\vdots&\ddots&0&0&0&\vdots\\
0&1&\ddots&\ddots&\vdots&\vdots&-a_1&-a_3&\ddots&\vdots&\vdots&\vdots&\vdots\\
\vdots&0&\ddots&\ddots&\vdots&\vdots&0&-a_1&\ddots&\vdots&\vdots&\vdots&\vdots\\
\vdots&\vdots&\ddots&\ddots&0&\vdots&0&0&\ddots&\ddots&\vdots&\vdots&\vdots\\
\vdots&\vdots&\ddots&\ddots&1&0&\vdots&\vdots&\ddots&0&-a_1&-a_3&\vdots\\
0&0&\cdots&\cdots&0&1&0&0&\cdots&0&0&-a_1&\times\\
\end{smallmatrix}\right)\end{equation*}
and here is the matrix of the multiplication by $x+y$
\begin{equation*}\left(\begin{smallmatrix}
0&0&\cdots&\cdots&0&\times&a_6&0&\cdots&0&\times&\times&\times\\
1&0&\ddots&\ddots&\vdots&\times&a_4&a_6&\cdots&0&\times&\times&\times\\
0&1&\ddots&\ddots&\vdots&\vdots&a_2&a_4&\ddots&\vdots&\vdots&\vdots&\vdots\\
\vdots&0&\ddots&\ddots&\vdots&\vdots&1&a_2&\ddots&\vdots&\vdots&\vdots&\vdots\\
\vdots&\vdots&\ddots&\ddots&\vdots&\vdots&0&\ddots&\ddots&\vdots&\vdots&\vdots&\vdots\\
0&\vdots&\ddots&\ddots&1&\times&\vdots&\ddots&\ddots&1&\times&\times&\vdots\\
1&0&\ddots&\ddots&0&0&-a_3&\vdots&\ddots&0&0&0&\vdots\\
0&1&\ddots&\ddots&\vdots&\vdots&1-a_1&-a_3&\ddots&\vdots&\vdots&\vdots&\vdots\\
\vdots&0&\ddots&\ddots&\vdots&\vdots&0&1-a_1&\ddots&\vdots&\vdots&\vdots&\vdots\\
\vdots&\vdots&\ddots&\ddots&0&\vdots&0&0&\ddots&\ddots&\vdots&\vdots&\vdots\\
\vdots&\vdots&\ddots&\ddots&1&0&\vdots&\vdots&\ddots&0&1-a_1&-a_3&\vdots\\
0&0&\cdots&\cdots&0&1&0&0&\cdots&0&0&1-a_1&\times\\
\end{smallmatrix}\right).\end{equation*}
\end{remark}

\begin{remark}
Proposition~\ref{minimum valuation} tells us that, if we take an equation for $E$ with coefficients in $\mathcal{O}_K$, and $u=ay+bx+c$ such that $a=1$ and $b,c\in\{0,1,-1\}$, then \[\chi_{u,\ell^k}\in \frac{1}{\ell^3}\mathcal{O}_K[X]\]
if $\ell$ is prime and $k$ a positive integer, and \[\chi_{u,n}\in \mathcal{O}_K[X]\] if $n$ is composite.
\end{remark}

\subsection{Case $n=3$}

Let $b_2=a_1^2+4a_2$, $b_4=2a_4+a_1a_3$, $b_6=a_3^2+4a_6$ and $b_8=a_1^2a_6+4a_2a_6-a_1a_3a_4+a_2a_3^2-a_4^2$.
Since a $3$-torsion point $(x,y)$ satisfies $[2](x,y)=-(x,y)$, it satisfies \[\frac{x^4-b_4x^2-2b_6x-b_8}{4x^3+b_2^2+2b_4x+b_6}=x.\]
Then, after simplification, \[\psi_3(x)=3x^4+b_2x^3+3b_4x^2+3b_6x+b_8.\]
If $a_1\neq0$, then $\chi_{y,3}$ has Galois group $\rho_{E,3}(G_K)$ from Corollary~\ref{pol réalisant im rho} and Subsection~\ref{condition on u}. We can compute explicitly $\chi_{y,3}$ by two ways: either as a characteristic polynomial, either as a resultant. In the first way, we compute the determinant of
\[\left(\begin{smallmatrix}-X&0&0&0&a_6&\frac{-b_8}{3}&\frac{-b_8}{3}(a_2-\frac{b_2}{3})&A_0\\
0&-X&0&0&a_4&a_6-b_6&\frac{-b_8}{3}-b_6(a_2-\frac{b_2}{3})&A_1\\
0&0&-X&0&a_2&a_4-b_4&a_6-b_6-b_4(a_2-\frac{b_2}{3})&A_2\\
0&0&0&-X&1&a_2-\frac{b_2}{3}&a_4-b_4-\frac{b_2}{3}(a_2-\frac{b_2}{3})&A_3\\
1&0&0&0&-X-a_3&0&0&\frac{a_1b_8}{3}\\
0&1&0&0&-a_1&-X-a_3&0&a_1b_6\\
0&0&1&0&0&-a_1&-X-a_3&a_1b_4\\
0&0&0&1&0&0&-a_1&-X+\frac{a_1b_2}{3}-a_3\end{smallmatrix}\right)\]
with \[A_0=-\frac{b_8}{3}\left(-\frac{b_2}{3}(a_2-\frac{b_2}{3})+a_4-b_4\right),\]
\[A_1=-b_6\left(-\frac{b_2}{3}(a_2-\frac{b_2}{3})+a_4-b_4\right)-\frac{b_8}{3}(a_2-\frac{b_2}{3}),\]
\[A_2=-b_4\left(-\frac{b_2}{3}(a_2-\frac{b_2}{3})+a_4-b_4\right)-b_6(a_2-\frac{b_2}{3})-b_8/3,\]
\[A_3=-\frac{b_2}{3}\left(-\frac{b_2}{3}(a_2-\frac{b_2}{3})+a_4-b_4\right)-b_4(a_2-\frac{b_2}{3})+a_6-b_6.\]

By the second way, we compute the determinant

\[\mathrm{Res}_x(\psi_3,w_E)=\det\left(\begin{smallmatrix}3&0&0&-1&0&0&0\\b_2&3&0&-a_2&-1&0&0\\3b_4&b_2&3&a_4-a_1y&-a_2&-1&0\\3b_6&3b_4&b_2&a_6-a_3-y^2&a_4-a_1y&-a_2&-1\\b_8&3b_6&3b_4&0&a_6-a_3-y^2&a_4-a_1y&a_2\\0&b_8&3b_6&0&0&a_6-a_3-y^2&a_4-a_1y\\0&0&b_8&0&0&0&a_6-a_3-y^2\end{smallmatrix}\right).\]

Both methods give the same polynomial $\chi_{y,3}$.
\begin{example}
Let $E:y^2+xy=x^3-\frac{4}{13}$. Its Galois image modulo $3$ is the Borel subgroup $\begin{pmatrix}*&*\\0&*\end{pmatrix}$ and is isomorphic to the Galois group of \[\chi_{y,3}(x)=x^8 - \frac{1}{3}x^7 - \frac{851}{351}x^6 + \frac{12}{13}x^5 + \frac{760}{507}x^4 + \frac{3076}{4563}x^3 - \frac{16}{169}x^2 + \frac{576}{2197}x - \frac{6912}{28561}\]
\end{example}

If $a_1=0$, we can take the function $u=x+y$ for example. With a short Weierstrass equation $E:y^2=x^2+Ax+B$, the matrix of the multiplication by $x+y$ is
\[\begin{pmatrix}
0&0&0&A^2/3&B&A^2/3&0&-\frac{A^3}{3}\\1&0&0&-4B&A&-3B&A^2/3&4BA\\0&1&0&-2A&0&-A&-3B&2A^2+A^2/3\\0&0&1&0&1&0&-A&-3B\\1&0&0&0&0&0&0&A^2/3\\0&1&0&0&1&0&0&-4B\\0&0&1&0&0&1&0&-2A\\0&0&0&1&0&0&1&0
\end{pmatrix}\]
and \[\mathrm{Res}_x(\psi_3',w_{E'})=\begin{pmatrix}
3&0&0&-1&0&0&0\\
0&3&0&-1&-1&0&0\\
6A&0&3&A-2y&-1&-1&0\\
12B&6A&0&B-y^2&A-2y&-1&-1\\
-A^2&12B&6A&0&B-y^2&A-2y&-1\\
0&-A^2&12B&0&0&B-y^2&A-2y\\
0&0&-A^2&0&0&0&B-y^2\end{pmatrix}.\]
We obtain \[\chi_{x+y,3}=x^8+(4A+8B)x^6+(\frac{-32}{3}A^2+8B)x^5+(\frac{8}{3}A^3+\frac{10}{3}A^2-40AB+18B^2)x^4\]
\[+(16AB-80B^2)x^3+(\frac{16}{3}A^4-\frac{4}{3}A^3+\frac{40}{3}A^2B+36AB^2+16B^2)x^2\]
\[+(\frac{-32}{9}A^4+\frac{32}{3}A^3B-\frac{8}{3}A^2B-16AB^2+72B^3)x\]
\[- \frac{16}{27}A^6 - \frac{8}{9}A^5 - 8A^3B^2 +\frac{1}{9} A^4 - \frac{8}{3}A^3B - 6A^2B^2 - 27B^4 - 16B^3).\]

By an appropriate change of variables, we can vary the valuation at a given prime as follows:

\begin{prop}Let $\got{p}$ be a prime ideal of $\mathcal{O}_K$ which does not contain $3$, and let $m$ be an integer. Let $E/K$ be an elliptic curve. We can choose $u$ such that, for $i=0,\dots,\deg(\chi_{u,3})$, the coefficient of degree $i$ of $\chi_{u,3}$ has valuation at $\got{p}$ greater or equal than $2m(\deg(\chi_{u,3})-i)$.
\end{prop}

\begin{proof}
Let $\pi\in\got{p}\setminus \got{p}^2$, and let $\lambda=\pi^m$. We have $v_\got{p}(\lambda)=m$. Let the equation $y^2=x^3+Ax+B$ be minimal for $E$ at $\got{p}$. Take $u=\lambda^3y+\lambda^2x$. Then $\chi_{u,3}$ is equal to $\chi_{x'+y',3}$ where $y'^2=x'^3+\lambda^4Ax'+\lambda^6B$. So\[\chi_{u,3}=x^8+4\lambda^4(A+2\lambda^2B)x^6+8\lambda^6\left(\frac{-4}{3}\lambda^2A^2+B\right)x^5\]\[+2\lambda^8\left(\frac{4}{3}\lambda^4A^3+\frac{5}{3}A^2-20\lambda^2AB+9\lambda^4B^2\right)x^4+16\lambda^{10}(AB-5\lambda^2B^2)x^3\]\[+4\lambda^{12}\left(\frac{4}{3}\lambda^4A^4-\frac{1}{3}A^3+\frac{10}{3}\lambda^2A^2B+9\lambda^4AB^2+4B^2\right)x^2 \]
\[+8\lambda^{14}\left(\frac{-4}{9}\lambda^2A^4+\frac{4}{3}\lambda^4A^3B-2\lambda^2AB^2+9\lambda^4AB^3-\frac{8}{3}A^2B\right)x\]
\scalebox{0.88}{$\begin{array}{r} +\lambda^{16}\left(\frac{-16}{27}\lambda^8A^6-\frac{8}{9}\lambda^4A^5-8\lambda^8A^4+\frac{1}{9}A^4-\frac{8}{3}\lambda^2A^3B-6\lambda^4A^2B^2-27\lambda^8B^4-16\lambda^2B^2\right).\end{array}$}
\end{proof}

The proof of the previous proposition suggests that we can obtain a polynomial with the smallest possible valuation at a prime $\got{p}$ taking a minimal equation at $\got{p}$.

\begin{remark}
Under the same notation as in the previous proposition, we can ask if a similar result holds for $n\geq5$, with some other equation linking the coefficients and the associated power of $\pi$. But, for that, we certainly have to compute $\chi_{x+y,n}$, for $n\geq5$ and we did not succeed with non-specialized coefficients for $E$. Nevertheless, we observe that, for $n=2$, a polynomial generating $K(E[2])$ is $\chi_{E,2}:=x^3+Ax+B$. If we make the change of variables $(A,B)\mapsto (\lambda^4A,\lambda^6B)$ with $\lambda=\pi^m$, we obtain $\chi_{E,2}=x^3+\lambda^4Ax+\lambda^6B$. This is an equation for $E$ such that, for all $i=0,\dots,\deg(\chi_{E,2})$, the coefficient of degree $i$ of $\chi_{E,2}$ is divisible by $\pi$ at least $2m(\deg(\chi_{E,2})-i)$ times.
\end{remark}

\section{More examples}\label{examples}

In the previous subsection, we have computed, using Sagemath \cite{sagemath}, the polynomial $\chi_{x+y,3}$, for any $E/K$ an elliptic curve with a short Weierstrass equation.

For $n$ larger than $3$, let us start by observing that $\deg\widetilde{\psi}_{5}=12=\deg\widetilde{\psi}_{6}$. So, thanks to Theorem~\ref{galois psi}, respectively Theorem~\ref{pol réalisant im rho}, we can construct polynomials of degree $12$, respectively degree $24$, with different Galois group. There are many others such cases, for example:
\[\deg\widetilde{\psi}_9=\deg\widetilde{\psi}_{10}=36,\]
\[\deg\widetilde{\psi}_{19}=\deg\widetilde{\psi}_{22}=180,\]
\[\deg\widetilde{\psi}_{31}=\deg\widetilde{\psi}_{33}=480,\]
\[\deg\widetilde{\psi}_{71}=\deg\widetilde{\psi}_{82}=2520.\]

The interesting point is that, thanks to these theorems, a Sagemath program is enough to compute such polynomials. We remark that an artefact of this method is to find the Galois group of a polynomial of high degree, just by calculating the characteristic polynomial of a matrix, which is easy (although sometimes long) for a computer, whereas numerically finding the Galois group of a polynomial of high degree is not feasible with the current technology. Thanks to Remark~\ref{two ways to calculate} it is even easier, since it suffices to compute a resultant of two well-known polynomials.

\begin{example}\label{ex serre curves}
\cite[Theorem 8.1]{Dserrecurves} gives us a curve $E/\Q(t)$ and a set $S$ such that the specialization $E_t$ of $E$ at $t$ is a Serre curve over $\Q$ if and only if $t\notin S$. In particular, from \cite[Theorem 1.8]{Dserrecurves}, for $t\notin S$, the representation $\rho_{E_t,n}$ is surjective for all prime, so for all product of pairwise distinct primes, as well as for $4$ and $9$. This curve is defined by\[E:y^2+xy=x^3+t.\]
Thanks to computations with Sagemath, we obtain, for each integer $n$ where $\rho_{E,n}$ is surjective (and for which the computation is feasible by a computer), a family of irreducible polynomials $(\chi_{u,n})_{t\notin S}$ with Galois group $\mathrm{GL}_2(\Z/n\Z)$. For example:
\[\scalebox{0.83}{$\begin{array}{r}\chi_{y_t,3}=x^8+\frac{1}{3}x^7 + (8t + \frac{1}{27})x^6 + 3tx^5 + (18t^2 + \frac{2}{3}t)x^4 + (-7t^2 + \frac{1}{27})x^3 - t^2x^2 + 9t^3x - 27t^4 \end{array}$}\]
has Galois group $\mathrm{GL}_2(\F_3)$ for all $t\notin S$, and
\[\scalebox{0.90}{$\begin{array}{r} \chi_{y_t,4}=x^{12} - \frac{1}{2} x^{11} + \left(54 t + \frac{1}{8}\right) x^{10} - \frac{55}{2} t x^{9} + \left(891 t^{2} + \frac{99}{8} t\right) x^{8} + \left(27 t^{2} - 2 t\right) x^{7}\\ + \left(2916 t^{3} - \frac{219}{4} t^{2} + \frac{1}{8} t\right) x^{6} + \left(-1863 t^{3} + 6 t^{2}\right) x^{5}\\ + \left(-24057 t^{4} + \frac{1107}{4} t^{3} + \frac{1}{8} t^{2}\right) x^{4} + \left(\frac{13851}{2} t^{4} + 4 t^{3}\right) x^{3}\\ + \left(39366 t^{5} - \frac{891}{8} t^{4} - \frac{1}{8} t^{3}\right) x^{2} - \frac{2187}{2} t^{5} x - 19683 t^{6} - \frac{729}{8} t^{5} - \frac{1}{8} t^{4}
\end{array}$}\]
has Galois group $\mathrm{GL}_2(\Z/4\Z)$.
\[\scalebox{0.80}{$\begin{array}{r} \chi_{y_t,5}=x^{24} - x^{23} + \left(216 t + \frac{3}{5}\right) x^{22} + \left(-217 t - \frac{3}{25}\right) x^{21} + \left(14742 t^{2} + \frac{877}{5} t + \frac{1}{125}\right) x^{20} \\+ \left(-7695 t^{2} - \frac{1971}{25} t\right) x^{19}+ \left(256608 t^{3} + 1477 t^{2} + \frac{506}{25} t\right) x^{18} + \left(-234495 t^{3} - \frac{3459}{25} t^{2} - 3 t\right) x^{17}\\+ \left(-\frac{19899513}{5} t^{4} + \frac{431487}{5} t^{3} + \frac{1287}{25} t^{2} + \frac{6}{25} t\right) x^{16}+ \left(3083670 t^{4} - \frac{324652}{25} t^{3} - \frac{448}{25} t^{2} - \frac{1}{125} t\right) x^{15}\\ + \left(\frac{79046928}{5} t^{5} - \frac{4052754}{5} t^{4} + 190 t^{3} + \frac{62}{25} t^{2}\right) x^{14} + \left(-\frac{23737698}{5} t^{5} + \frac{2105946}{25} t^{4} + \frac{3531}{25} t^{3} - \frac{3}{25} t^{2}\right) x^{13}\\ + \left(-\frac{88258572}{5} t^{6} - \frac{3661038}{5} t^{5} - \frac{15616}{5} t^{4} - \frac{252}{25} t^{3}\right) x^{12}\\ + \left(-\frac{36256086}{5} t^{6} + \frac{10623798}{25} t^{5} + \frac{11674}{25} t^{4} - \frac{1}{25} t^{3}\right) x^{11}\\+ \left(-17006112 t^{7} + \frac{29406402}{5} t^{6} - \frac{6267294}{125} t^{5} - \frac{11448}{125} t^{4}\right) x^{10}\\ + \left(48931938 t^{7} - \frac{22323438}{25} t^{6} + \frac{41099}{25} t^{5} + 5 t^{4}\right) x^{9}\\+ \left(\frac{798755823}{25} t^{8} - \frac{384172794}{25} t^{7} + \frac{216756}{25} t^{6} - \frac{1299}{25} t^{5} - \frac{3}{25} t^{4}\right) x^{8}\\ + \left(-\frac{463947993}{5} t^{8} + \frac{26690148}{25} t^{7} + \frac{103356}{25} t^{6} + \frac{36}{5} t^{5}\right) x^{7}\\+ \left(\frac{114791256}{25} t^{9} + \frac{457452603}{25} t^{8} + \frac{1010394}{25} t^{7} - \frac{4106}{25} t^{6} - \frac{13}{25} t^{5}\right) x^{6}\\+ \left(\frac{2157119019}{25} t^{9} - \frac{80326323}{125} t^{8} - \frac{796311}{125} t^{7} - \frac{1046}{125} t^{6} + \frac{1}{125} t^{5}\right) x^{5}\\+ \left(-\frac{86093442}{5} t^{10} - \frac{320458923}{25} t^{9} - \frac{1712421}{25} t^{8} + \frac{54}{25} t^{7} + \frac{6}{25} t^{6}\right) x^{4}\\+ \left(-\frac{186535791}{5} t^{10} + \frac{6121413}{25} t^{9} + \frac{13122}{5} t^{8} + \frac{22}{5} t^{7}\right) x^{3}+ \left(\frac{81310473}{25} t^{10} + \frac{236196}{25} t^{9} - \frac{324}{25} t^{8} - \frac{1}{25} t^{7}\right) x^{2}\\+ \left(\frac{43046721}{5} t^{11} + \frac{531441}{25} t^{10}\right) x+ \frac{387420489}{125} t^{12} + \frac{14348907}{125} t^{11} - \frac{19683}{25} t^{10} - \frac{729}{125} t^{9} - \frac{1}{125} t^{8}
\end{array}$}\]
has Galois group $\mathrm{GL}_2(\F_5)$, whereas
\[\scalebox{0.80}{$\begin{array}{r} \chi_{y_t,6}=x^{24} - x^{23} + \left(648 t + 1\right) x^{22} - 649 t x^{21} + \left(132678 t^{2} + 875 t\right) x^{20}+ \left(-68607 t^{2} - 462 t\right) x^{19}\\ + \left(6940080 t^{3} + 25075 t^{2} + 136 t\right) x^{18} + \left(-6152031 t^{3} - 8067 t^{2} - 19 t\right) x^{17}\\  + \left(-317375253 t^{4} + 1723113 t^{3} + 2832 t^{2} + t\right) x^{16}+ \left(256565718 t^{4} + 157736 t^{3} - 400 t^{2}\right) x^{15}\\  + \left(4311049392 t^{5} - 57426246 t^{4} - 117566 t^{3} - 4 t^{2}\right) x^{14}\\ + \left(-2135382426 t^{5} + 2165940 t^{4} + 12647 t^{3} + 3 t^{2}\right) x^{13}\\ + \left(-26822181564 t^{6} + 246120606 t^{5} + 413098 t^{4} - 152 t^{3}\right) x^{12}\\ + \left(10811281410 t^{6} + 22940172 t^{5} + 1918 t^{4} - 33 t^{3}\right) x^{11}\\ + \left(85506731136 t^{7} - 1285601706 t^{6} - 7285302 t^{5} - 9201 t^{4} + 3 t^{3}\right) x^{10}\\ + \left(-30700637982 t^{7} - 33424650 t^{6} + 561395 t^{5} + 1198 t^{4}\right) x^{9}\\ + \left(-159947266329 t^{8} + 3317254722 t^{7} + 20098530 t^{6} + 7890 t^{5} - 49 t^{4}\right) x^{8}\\ + \left(53354019195 t^{8} + 71226216 t^{7} - 1074060 t^{6} - 1875 t^{5} + t^{4}\right) x^{7}\\ + \left(185847043464 t^{9} - 5354268075 t^{8} - 26764506 t^{7} + 4618 t^{6} + 86 t^{5}\right) x^{6}\\ + \left(-46476109773 t^{9} + 172186884 t^{8} + 1826145 t^{7} + 2342 t^{6} - t^{5}\right) x^{5}\\ + \left(-134047489194 t^{10} + 2903262183 t^{9} + 3136158 t^{8} - 47952 t^{7} - 89 t^{6}\right) x^{4}\\ + \left(12038732973 t^{10} - 115145550 t^{9} - 572994 t^{8} - 136 t^{7} + t^{6}\right) x^{3}\\ + \left(55788550416 t^{11} - 1076168025 t^{10} - 2480058 t^{9} + 8019 t^{8} + 17 t^{7}\right) x^{2}\\ + \left(1937102445 t^{11} + 87687765 t^{10} + 367416 t^{9} + 405 t^{8}\right) x\\- 10460353203 t^{12} - 243931419 t^{11} - 2480058 t^{10} - 8019 t^{9} - 8 t^{8}
\end{array}$}\]
has Galois group $\mathrm{GL}_2(\Z/6\Z)$. As we have underlined before, the polynomials $\chi_{y_t,5}$ and $\chi_{y_t,6}$ both have degree $24$. 
\end{example}

The following two examples were found using the database \cite{lmfdb} and come from \cite{sutherland_2016}.

\begin{example}
The elliptic curve \[E:y^2+xy=x^3-x^2-9x+3699\]
is defined over $\Q$, and has surjective Galois image over $\Q$ for all primes except for $7$. We have \[G:=\rho_{E,7}(G_\Q)=\left\langle\begin{pmatrix}6&0\\0&5\end{pmatrix},\begin{pmatrix}6&6\\0&4\end{pmatrix}\right\rangle.\]Over $K=\Q(\sqrt{-3})$, the Galois representation is also surjective for all primes except for $7$ and we have \[G':=\rho_{E,7}(G_K)=\left\langle\begin{pmatrix}1&6\\0&2\end{pmatrix},\begin{pmatrix}1&6\\0&5\end{pmatrix}\right\rangle.\]
Hence, $\chi_{x+y,7}$, which has degree $48$, has Galois group $G$, which has order $84$, over $\Q$ and Galois group $G'$, which has order $42$, over $\Q(\sqrt{-3})$. In particular, $\chi_{x+y,7}$ has a rational root over $\Q(\sqrt{-3})$. Since $G'$ does not contains $-\id$, the polynomial $\psi_7$, which has degree $24$, also has Galois group $G$ over $K$. It also has a root in $\Q(\sqrt{-3})$.
\end{example}

\begin{example}
Let be the elliptic curve $E:y^2=f(x)$ where \[f(x)=x^3+ix^2+(2i-2)x-2i-1.\] It is defined over $K=\Q(i)$. The image of $\rho_{E,\ell}$ is surjective for all primes $\ell$ except $2$ and $5$. The image mod $5$ is $\left\langle\begin{pmatrix}0&3\\4&0\end{pmatrix},\begin{pmatrix}3&3\\2&3\end{pmatrix}\right\rangle$. Hence, for all primes $\ell$ except $2$ and $5$, the polynomial $\chi_{x+y,\ell}$, which has degree $\ell^2-1$, is irreducible and has Galois group $\mathrm{GL}_2(\F_\ell)$, which has order $(\ell^2-1)(\ell-1)\ell$. The polynomial $\chi_{x+y,5}$, which has degree $24$, has Galois group $\left\langle\begin{pmatrix}0&3\\4&0\end{pmatrix},\begin{pmatrix}3&3\\2&3\end{pmatrix}\right\rangle$, which has order $96$. The image of $\rho_{E,2}$ is $\left\langle\begin{pmatrix}1&1\\1&0\end{pmatrix}\right\rangle$, which has order 3, and is the Galois group of $f(x)$. %In particular, $f(x)$ has a rational roots over $\Q(i)$.
We can find the image modulo $10$ from the image modulo $2$ and modulo $5$. Using that the reduction modulo $5$, respectively modulo $2$, of $\rho_{E,10}(G_K)$ is $\rho_{E,5}(G_K)$, respectively $\rho_{E,2}(G_K)$, we find that \[\rho_{E,10}(G_K)=\left\langle \begin{pmatrix}3&3\\7&8\end{pmatrix},\begin{pmatrix} 0&3\\9&5\end{pmatrix}\right\rangle.\]Hence $\chi_{x+y,10}$, which has degree $72$, has Galois group $\left\langle \begin{pmatrix}3&3\\7&8\end{pmatrix},\begin{pmatrix} 0&3\\9&5\end{pmatrix}\right\rangle$ which has order $288$.
\end{example}

\bibliographystyle{annotate}
\bibliography{Ref}

\end{document}